\documentclass{gtart} 

\input gtoutput
\volumenumber{4}\papernumber{2}\volumeyear{2000}
\pagenumbers{85}{116}\published{28 January 2000}
\proposed{Walter Neumann}\seconded{Joan Birman, Wolfgang Metzler}
\received{9 February 1999}\revised{10 November 1999}
\accepted{13 January 1999}

\usepackage{amssymb}
\let\Bbb\mathbf

\def\S{Section }

\newcounter{punkt}[section]
\newcommand{\pu}{\ppar\refstepcounter{punkt}%
{\bf\arabic{section}.\arabic{punkt}}\qua}

\def\bokvnack{{\frac{\bokv}{c_k}}}

\def\oneoverlambda{{\frac{1}{\lambda}}}
\def\onenacj{{\frac{1}{c_j}}}

\def\onenack{{\frac{1}{c_k}}}

\def\calcspec{{\calc_{spec}}}

\def\deltaonespec{{\Delta^{(1)}_{spec}}}
\def\deltaone{{\Delta^{(1)}}}
\def\deltazerospec{{\Delta^{(0)}_{spec}}}
\def\coxsph{{\Sigma_{sph}}}

\def\widewa{{\widehat{W_a}}}

\def\chxy{{\mbox{ch}(x,y)}}

\def\eonen{e_1 ,\ldots,e_n}
\def\xon{x ,x_1,\ldots,x_n}

\def\R{{\Bbb{R}}}
\def\Rp{{\R_{\geq 0}}}
\def\Z{{\Bbb{Z}}}

\def\N{{\Bbb{N}}}

\def\Rnone{{\R^{n+1}}}

\def\Rn{{\R^n}}

\def\cala{{\cal A}}

\def\calc{{\cal C}}
\def\calg{{\cal G}}
\def\calh{{\cal H}}
\def\cA{{\cal A}}

\def\eps{{\epsilon}}
\def\epsone{{\eps_1}}
\def\epszero{{\eps_0}}
\def\epsn{{\eps_n}}
\def\epsi{{\eps_i}}
\def\epsj{{\eps_j}}
\def\epsk{{\eps_k}}
\def\epsir{{\eps_{i_r}}}
\def\epsij{{\eps_{i_j}}}
\def\epsik{{\eps_{i_k}}}

\def\lemma{{\bf Lemma}\qua}
\def\proof{\par\medskip{\bf Proof}\qua}

\def\definitions{{\bf Definitions}\qua}
\def\theorem{{\bf Theorem}\qua}

\def\remarks{{\bf Remarks}\qua}

\def\corollary{{\bf Corollary}\qua}

\def\Phivee{{\Phi^\vee}}
\def\alphavee{{\alpha^\vee}}

\def\ZPhivee{{\Z\Phi^\vee}}
\def\ZPhid{{\Z\Phi^\#}}

\def\alphavee{{\alpha^\vee}}

\def\alphai{{\alpha_i}}

\def\boi{{\bar{\omega}_i}}

\def\boov{{\bar{\omega}_1^\vee}}
\def\botv{{\bar{\omega}_2^\vee}}
\def\bonv{{\bar{\omega}_n^\vee}}
\def\bojv{{\bar{\omega}_j^\vee}}
\def\bokv{{\bar{\omega}_k^\vee}}

\def\bov{{\bar{\omega}^\vee}}

\def\boiov{{\bar{\omega}_{i_1}^\vee}}

\def\boirv{{\bar{\omega}_{i_r}^\vee}}

\def\boiv{{\bar{\omega}_i^\vee}}

\begin{document}
\author{Gennady A Noskov}
\title{Combing Euclidean buildings}

\address{Russia, 644099, Omsk\\Pevtsova 13, IITAM SORAN}
\email{noskov@private.omsk.su}

\begin{abstract}
For an arbitrary Euclidean building   we define a certain
combing, which satisfies the ``fellow traveller property''
and admits a recursive definition.
 Using this combing we prove that any group acting freely, cocompactly
and by order preserving automorphisms
on a  Euclidean building of one of the types $A_n,B_n,C_n$ 
admits a  biautomatic structure.
\end{abstract}

\asciiabstract{For an arbitrary Euclidean building we define a certain
combing, which satisfies the ``fellow traveller property'' and admits
a recursive definition.  Using this combing we prove that any group
acting freely, cocompactly and by order preserving automorphisms on a
Euclidean building of one of the types A_n,B_n,C_n admits a
biautomatic structure.}

\primaryclass{20F32}

\secondaryclass{20F10}

\keywords{Euclidean building, automatic group, combing}

\maketitlepage

\section{Introduction}

 Let $G$ be a group that acts properly and cocompactly on a  
piecewise Euclidean simply connected $CAT(0)$--complex $\Delta$
(see eg \cite{bridsontrieste} for definitions).
 (The action of course is supposed to be cellular, properness
means that the isotropy group $G_\sigma$ is finite for every
cell $\sigma$ and cocompactness means that $\Delta$ has only finitely
many cells ${\rm mod}~G$.)
  It is still unknown whether $G$ is (bi)automatic.
 Moreover,  the question remains unanswered even  in the case when
$\Delta$ is a Euclidean building \cite{cartwrightshapiro}.
 
 ``It is reasonable to guess that the answer is `yes' because of
the work of Gersten and Short and because of the geometry and regularity
present in buildings, but this is far from a trivial question''
(John Meier's review [MR 96k:20071] of the paper
\cite{cartwrightshapiro}).

 The first results in this direction
are  contained in the papers of S Gersten and H Short 
\cite{gerstenshort1}, \cite{gerstenshort2}
where it is proven that if 
$G$ is given by a finite presentation
satisfying the small cancellation conditions $C(p),T(q)$ 
$((p,q)=(6,3),(4,4),(3,6))$ then $G$ is biautomatic.
  They showed in \cite{gerstenshort1} that the fundamental 
group of a piecewise Euclidean $2$--complex of nonpositive curvature of type 
$A\sb 1\times A\sb 1$ or $A\sb 2$ is automatic. 
 ($A\sb 1\times A\sb 1$ corresponds to the Euclidean planar 
tessellation by
unit squares, and $A\sb 2$ to the tessellation by equilateral triangles). 
 In the subsequent paper \cite{gerstenshort2} the
authors prove an analogous result for $2$--complexes of types $B\sb 2$ and 
$G\sb 2$ corresponding to the Euclidean tessellations by 
$(\frac \pi 2, \frac \pi 4,\frac \pi 4)$ and $(\frac \pi
2,\frac \pi 3,\frac \pi 6)$ triangles, respectively.
  It follows from this work
that any torsion free group $G$ which admits a proper  cocompact action
on a  Euclidean building of type $A\sb 2$ is biautomatic.
 W Ballmann  and M Brin \cite{ballmannbrinpolygonal} 
have proven the automatic property   for a group $G$ which acts
simply transitively on the vertices of a  simply connected
 (3,6)--complex.
 D Cartwright and M Shapiro have proven the following theorem 
\cite{cartwrightshapiro}:
  Let $G$ act simply transitively on the vertices of a
Euclidean  $A\sb n$--building in a type rotating way.
 Then $G$ admits a
geodesic, symmetric automatic structure.
 In \cite{mytriangle} some variation of this result is proven in a
more geometric way in the case of $n=2$.  
  It  is worth mentioning that in the case of nonpositively curved cube 
complexes the general result was obtained by G Niblo and L Reeves 
\cite{nibloreevescubeaut}, 
namely,  any group acting properly and cocompactly on such a complex is 
biautomatic.

 In this paper we  define a certain 
combing on an arbitrary Euclidean building, prove the
``fellow traveller property'' for this combing and 
the ``recursiveness property''.
Our main result is the following.\ppar\goodbreak

\theorem\par\nobreak
{\sl 
{\rm (1)}\qua Let $\Delta$ be any  Euclidean building of one of the
types $A_n,B_n,C_n$, ordered in a standard way  
(see \S\ref{standard ordering} for a definition).
 Then any group acting   freely and cocompactly
on $\Delta$ by order-preserving automorphisms
admits a  biautomatic structure.

{\rm (2)}\qua
 If $\Delta$ is any  Euclidean building of one of the
types $A_n,B_n,C_n$, then any group acting   freely and 
cocompactly on $\Delta$ is virtually biautomatic (that is there
is a finite index subgroup in it, possessing a biautomatic structure).
}\medskip 

 In \S\ref{euclidean coxeter complexes} we review some of the standard 
facts on Euclidean Coxeter complexes.
 In \S\ref{ordering buildings} we introduce  the main notion of an ordering
of a Euclidean building and prove that any Euclidean building can be ordered.
 In \S\ref{sec def combing} we define a natural combing $\calc$ on a 
Euclidean building .
 In Sections \ref{ftp} and \ref{recursive} we prove the
``fellow traveller property''
and  the ``recursiveness property'' for a combing $\calc$.
 The concluding \S\ref{main} is devoted to the proof our main result.

\medskip
{\bf Acknowledgements}\qua
I am grateful to the SFB 343, University of Bielefeld for their
hospitality in the falls of 1996--97 years while I carried out most  of this 
work.
 I would like to thank
Herbert Abels for his kind invitation, interest and support.
 Thanks to Sarah Rees for the help to  make the language
of the paper  more regular.
 Thanks to referee for many improvements of the text.
 The work was supported in part by a RFFI grant N 96-01-01610.

\section{Euclidean Coxeter complexes}
\label{euclidean coxeter complexes}
 For the convenience of the reader we recall the relevant material from
\cite{bourbaki456}, \cite{humphreysreflection}, thus making our exposition 
self-contained.

\pu {\bf Roots and Weyl group}\par\label{roots}
 Let $\Phi$ to be a {\it root system}, which is supposed to be
reduced, irreducible and crystallographic.
 That is $\Phi$ is a finite set of nonzero  vectors, spanning a finite
dimensional Euclidean space $V$ and such that

1)\qua$\Phi \cap  \R \alpha = \{\alpha , -\alpha\}$ for all
       $\alpha \in \Phi$,

2)\qua $\Phi$ is invariant under reflection $s_{\alpha}$ in the hyperplane
$ H_{\alpha}$ orthogonal to $\alpha$ for all $\alpha \in \Phi$,

3)\qua ${\frac {2(\alpha,\beta)}{(\alpha,\alpha)}} \in \Z$ for all
   $\alpha,\beta \in \Phi$,

4)\qua $V$ does not  admit an orthogonal
decomposition $V=V'\oplus V''$ such that  $\Phi=\Phi'\cup\Phi''$ with
$\Phi\subset V',\Phi''\subset V''$.

 The {\it Weyl group} $W$ of $\Phi$ is the group generated by
all reflections $s_{\alpha}~(\alpha \in \Phi)$.
 In equal terms $W$ is generated by all reflections $s_H$, where $H$ ranges
over the set ${\cal H}$ of all hyperplanes, orthogonal to the
roots from $\Phi$.

For any choice of the basis $\eonen$ of $V$ there is a
lexicographic order on $V$, where
$\sum a_{i}e_{i} < \sum b_{i}e_{i}$ means that $a_{k} < b_{k} $ if $k$ is the
least index $i$ for which $a_{i} \not= b_{i} $.
 We call a subset $\Phi_+ \subset \Phi $ a {\it positive system} if it  consists
of all those roots which are positive relative to some ordering of $V$ of the kind
above.

 If this is the case, then $\Phi$ must be the disjoint union of $\Phi_+$ and
$-\Phi_+$, the latter being called a {\it negative system}.
 When $\Phi_+$ is fixed, we can write
$\alpha> 0$ in place of $\alpha \in \Phi_+$.
 It is clear that positive systems exist.

 Call a subset $\Pi$ of $\Phi$ a {\it simple system} if $\Pi$ is a
vector space basis for  $V$ and if moreover each $\alpha \in \Phi$ is a linear
combination of $\Pi$ with coefficients all of the same sign (all
nonnegative or all nonpositive).
 If $\Pi$ is a simple system in $\Phi$, then there is a unique positive system
containing $\Pi$.
 Every positive system $\Pi$ in $\Phi$  contains a unique simple system; in
particular, simple systems exist.

 Any two positive (resp.\ simple) systems in $\Phi$ are conjugate under $W$.
 Thus $W$ permutes the the various positive (or simple) systems in a transitive
fashion.
 This permutation action is indeed a simply transitive action, that is if
$w\in W$ leaves  the positive (or simple) system invariant, then $w=1.$

\pu{\bf Coroots, lattices}\par
 Setting $\alphavee:=2\alpha/(\alpha,\alpha)$, the set $\Phivee$ of
all {\it coroots} $\alphavee(\alpha \in \Phi)$ is also a root system in
$V$, with simple system $\Pi^{\vee}:=\{\alphavee|\alpha \in \Pi \}$.
 The Weyl group of $\Phivee$ is $W$, with $w\alphavee=(w\alpha)^{\vee}$.
 The $\Z$--span $\Z\Phi$ of $\Phi$ in $V$ is called the {\it root lattice}; it
is a lattice in $V$.
 Similarly, we define the {\it coroot lattice} $\Z\Phi^\vee$.
 Define the {\it coweight lattice} $\ZPhid$~-- it is just a dual lattice of the
root  lattice $\Z\Phi$, that is
$$
\ZPhid =  \{\lambda\in V|(\lambda,\alpha)\in\Z~\mbox{for all}~
\alpha\in\Phi \}.
$$
 Since $(\Phivee,\Phi)\subset \Z$ and both $\ZPhivee,\Z\Phi$ are the lattices,
one can conclude that $\ZPhid$ contains $\ZPhivee$ as a subgroup of finite
index.

\pu{\bf Fundamental domain and spherical Coxeter complex}\par
\label{spherical cox}
 Let $W$ be the Weyl group of a root system $\Phi$.
 The hyperplanes $H$ with $s_H \in W$ cut $V$ into polyhedral pieces, which
turn out to be cones over simplices.
 One obtains in this way a simplicial complex 
$\coxsph=\coxsph(W)$ which triangulates 
the
unit sphere in $V$.
 This is a {\it spherical Coxeter complex}.
 More exactly let $\Phi_+$ be  a positive system, containing the simple system
$\Pi$.
 Associated with each hyperplane $H_{\alpha}$ are the closed
half-spaces $H_{\alpha}^+$ and $H_{\alpha}^-$, where
$H_{\alpha}^{+} = \{\lambda \in V|(\lambda,\alpha)\geq 0 \}$ and
$H_{\alpha}^{-} = \{\lambda \in V|(\lambda,\alpha)\leq 0 \}$.
 Define a {\it sector} $S=S_\Pi := \cap_{\alpha \in \Pi}H_{\alpha}^+ $
{\it associated to} $\Pi$.
 As an intersection of closed  convex subsets, $S$ is itself closed
 and convex.
 It is also a cone (closed under nonnegative scalar multiples).
 Sectors associated to $W$ are always {\it simplicial cones}, by which we
mean that, for some basis $\eonen$ of $V$ the sector $S$ consists
of the linear combinations $\sum a_{i}e_{i}$ with all $a_{i}$ positive.
 (In other words, $S$ is a cone
over the closed simplex with vertices $\eonen$).
 We call $\Rp e_i$  the {\it defining rays} of $S$.
 One can describe the defining rays of the sector $S_\Pi$  more explicitly in
terms of the basis of coroot lattice.
 Namely, let $\{\boiv \}$ be the dual basis of
$\Pi=\{\alpha_1,\ldots,\alpha_n\}$, that is
$(\boiv,\alpha_j)=\delta_{ij}$ for all $i,j=1,\ldots,n$.
 Then
$$
\sum_{1\leq i\leq n}a_i\boiv\in S\Longleftrightarrow
 a_j=(\sum_{1\leq i\leq n}a_i\boiv,\alpha_j)\geq 0,~j=1,\ldots,n
$$
We assert that the rays  $\Rp\boiv$ are the  defining rays for $S_\Pi.$
 Indeed, each line $\R~\boiv$ is precisely the line obtained by intersecting
all but one $H_{\alpha_i}$, namely
$
\R~\boiv=\cap_{j\not=i}H_{\alpha_j}.
$
  Consequently one of the halflines of this line is a defining ray and calculating
the scalar products we conclude that this is exactly $\Rp\boiv$.


 $W$ acts {\it simply transitively} on simple systems and this translates 
into a
simply transitive action on the the  sectors.
 This means that any two sectors are conjugate under
the action of $W$ and if $wS=S$ then $w=1$.
 Moreover any sector $S$ is a fundamental domain of the action
of $W$ on $V$, ie, each $\lambda\in V$ is conjugated under $W$ to one and 
only one point in $S$.
 The sectors are characterized topologically as the closure of the 
connected components of
the complement in $V$ of $\cup H_{\alpha}$ .
 They are in one one correspondence with the top-dimensional simplices
(= {\it chambers})
of the corresponding spherical complex.
 Given a sector $S_\Pi$ corresponding to a simple
system $\Pi $, its {\it walls} are defined to be the hyperplanes
$H_{\alpha}~(\alpha \in \Pi)$ .

\pu{\bf  Euclidean reflections and Euclidean Weyl group}\par
Let $\Phi$ be the root system in $V$ as it was defined in \S\ref{roots}.
 For each root $\alpha$ and each integer $k$, define a Euclidean hyperplane
$
H_{\alpha,k}: = \{\lambda \in V|(\lambda,\alpha) = k \}.
$
 Note that $H_{\alpha,k}=H_{-\alpha,-k} $ and that $H_{\alpha,0}$ coincides
with the reflecting hyperplane $H_{\alpha}.$
 Note too that $H_{\alpha,k}$ can be obtained by translating $H_{\alpha}$ by
$\frac{k}{2}\alphavee$.
 Define the corresponding Euclidean reflection as follows:
$s_{\alpha,k}(\lambda)
:= \lambda- ((\lambda,\alpha) - k)\alphavee.
$
We can also write $s_{\alpha,k}$ as $t(k\alphavee)s_{\alpha}$,
where $t(\lambda)$ denotes the translation by a vector $\lambda$.
 In particular, $s_{\alpha,0}=s_\alpha$.
 Denote by $\calh_\Z$ the collection of all hyperplanes
$H_{\alpha,k}~(\alpha \in \Phi,k \in \Z )$ which we shall call
the {\it walls}.
The elements of $\calh_\Z$ are permuted
in a natural way by $W_a$
as well as by  translations $t(\lambda)$, where $\lambda\in V$
satisfies $(\lambda,\mu)\in \Z$ 
for all roots $\alpha$ (that is $\lambda\in\ZPhid$).
In particular, $\ZPhid$
permutes the hyperplanes in $\calh_\Z$, hence so does its subgroup $\ZPhivee$.
 Define the {\it affine Weyl group} $W_a$ to be the subgroup of 
{\it Aff}$(V)$ generated by all affine reflections $s_{\alpha,k}$ where 
$a\in \Phi, k\in \Z$.
 Another description of $W_a$ is that it is the semidirect product
$W_a=\Z\Phi^\vee \rtimes W$ of the
finite Weyl group $W$ and the translation group
corresponding to the coroot lattice  $\Z\Phi^\vee$, see 
\cite{humphreysreflection}, \S 4.2.

 Since the translation group corresponding to $\ZPhid$ is also normalized
by $W$, we can form the semidirect product
$\widehat{W_a}=\ZPhid \rtimes W$, which
contains $W_a$ as a normal subgroup of finite index.
 Indeed, $\widehat{W_a}/W_a$ is isomorphic to $\ZPhid/\Z\Phivee$.
 One can easily see from 1), 2) that $\widehat{W_a}$ also permutes the 
hyperplanes in $\calh_\Z$.
 We call this group the {\it extended affine Weyl group}.

 \pu{\bf Euclidean Coxeter complexes}\par
 The hyperplanes $H\in\calh_\Z$ triangulate the space $V$ and the resulting
piecewise Euclidean complex
$\Sigma=\Sigma_\Phi$
is a {\em Euclidean Coxeter complex.}
 More generally we shall apply the same term to the Euclidean
simplicial structure $\Sigma$ on a Euclidean space $V'$ such that
that for some root system $\Phi$ in a Euclidean space $V$ there
is an affine isometry $\phi:V\rightarrow V'$ which induces simplicial
isomorphism between $\Sigma_\Phi$ and $\Sigma$.
 In particular in $\Sigma$ we have all the notions as in $\Sigma_\Phi$.
The extended Weyl group $\widewa$ acts by simplicial isometries on
$\Sigma_\Phi$ and this translates by $\phi$ to the action on $\Sigma$
but not in a  canonical way --  if $\phi:V\rightarrow V'$ is another 
isometry then the actions are conjugate by a suitable  isometry
of $V$.
 The possible ambiguity is resolved by the following lemma.

 \pu
\label{weyl groups are invariant}
\lemma
{\sl  Both $\widewa$ and $W_a$ are invariant under the
conjugation by any isometry $\phi$ of $V$, which preserves the
simplicial structure $\Sigma_\Phi$.}\medskip

 In particular the images of $\widewa$ and $W_a$ in ${\rm Aut}(\Sigma)$ 
are canonically defined and we call them the {\it extended affine Weyl  
group of $\Sigma$} and by the {\it affine Weyl  group of $\Sigma$} 
respectively.

\proof
 Since  $\phi$ leaves invariant the family of hyperplanes
$\calh_\Z$, it also   leaves invariant the family of reflections 
in the hyperplanes
of this family, hence normalizes the affine Weyl group $W_a$.
 Next it leaves invariant the set of special vertices 
(see definition in \S\ref{def special} and lemma \ref{special are zphid}) 
hence normalizes the translation group $\ZPhid$.
 Since $\widewa$ is generated by $W_a$ and $\ZPhid$ it is also 
normalized by $\phi$.
\qed

The collection $\cA$ of top-dimensional closed simplices consists of
the closures of the connected components of 
$V^{\circ}:= V\setminus \cup_{H\in\calh_\Z}H$.
 Each element of $\cA$ is called an {\it  alcove}.
 The group $W_a$ acts simply transitively on $\cA$,
\cite{humphreysreflection}, Chapter 4, Theorem 4.5.
 Any alcove $A$ is a fundamental domain of the action
of $W_a$ on $V$, ie, each $\lambda\in V$ is conjugated under $W$ to one and
only one point in $A$.
 In particular $V=\cup\{S_\Pi : \Pi~ {\rm is~ a~ simple~ system}\}.$
 Since $\widehat{W_a}$ permutes the hyperplanes
in $\calh_\Z$, it acts  simplicially on $\Sigma$.

\pu
{\bf Standard alcove}\par
\label{nice alcove}
 There is an alcove with a particularly nice description
(see
\cite{bourbaki456}, Corollary of Proposition 4 in \S 2, Chapter VI or
\cite{humphreysreflection} \S 4.9 ).
 Namely let $\Pi=\{\alphai\}$ be a  simple root system for $\Phi$.
 Let $\{\boiv\}$ be the dual basis for $\{\alpha_i\}$ -- this is the 
basis of the coweight lattice $\ZPhid$.
 Let $\tilde{\alpha}=\sum_{1\leq i \leq n}c_i\alpha_i$ be the 
corresponding highest root.
 Then the alcove $A=A_\Pi$, associated to $\Pi$ is a closed simplex with the
vertices $0$ and $\frac{1}{c_i}\boiv, i=1,\ldots,n$.
 We call this alcove  {\em a standard alcove associated to $\Pi$}.
 
 Another description of $A=A_\Pi$ is given by the formula
$
A= \cap_{\alpha \in \Pi}H_{\alpha}^+ \cap H^-_{\tilde{\alpha},1}, 
$
where $H^-_{\tilde{\alpha},1}$ is a closed negative half-space defined 
by the hyperplane $H_{\tilde{\alpha},1}$ and 
$\tilde{\alpha}$
is the highest root.
 Comparing this description with the definition of a sector 
$S_\Pi$ given in \S \ref{spherical cox}
we found that $A_\Pi$ sits on the top of the sector $S_\Pi$ and the defining
rays of $S_\Pi$ correspond to the ordered edges of $\Sigma$ having 0 as 
its origin.
 There is a one one correspondence  $A_\Pi\leftrightarrow S_\Pi$ between
the set $\cala_0$ of alcoves having 0 as a vertex and the set of sectors 
of a spherical Coxeter complex $\coxsph$.
 $W$ acts on $\cala_0$ and the fact that any $S_\Pi$ is a fundamental
domain for the action of $W$ on $\coxsph$ translates to the fact that
$A_\Pi$ is a fundamental domain for the action of $W$ on $\cala_0$ in a sense
that any directed edge of $A_\Pi$ having 0 as its origin is $W$--conjugate to 
one and only one such a vertex of $A_\Pi$.  

 \pu\label{def special}
{\bf Special vertices}\par
 The vertex $x\in\Sigma_\Phi$ is called a {\it special vertex} if its 
stabilizer $S_{W_a}(x)$ in $W_a$ maps isomorphically onto the associated finite
Weyl group $W$.
 (Note that the stabilizers of any vertex in $W_a$ and in $\widehat{W_a}$
coincide).
 Equivalently, for any hyperplane $H\in\calh_\Z$ there is a parallel 
hyperplane in $\calh_\Z$, passing through $x$.
 Yet another equivalent definition is that the maximal possible number
of hyperplanes from $\calh_\Z$ pass through $x.$

 \pu
\label{special}
\lemma\label{special are zphid}
{\sl
 The  set of special vertices  of the complex $\Sigma$ coincide with
the lattice  $\ZPhid$ (see
\cite{bourbaki456}, Proposition 3 in \S 2, Chapter VI).
}
\proof
 Since the zero vertex is special and  the coweight lattice 
$\ZPhid$ acts simplicially on $\Sigma$,
we conclude that $\ZPhid$ consists of the special vertices.
 Conversely, let $x$ be a special vertex.
 Since $W_a$ preserves the property of the vertex being special
and since it acts transitively on the set of alcoves,
we may assume that $x$ is the vertex of  the standard alcove
$$
A=<0,\frac{1}{c_1}\boov,\ldots,\frac{1}{c_n}\bonv>
$$
described above.
 If $x=0$, then obviously $x\in\ZPhid$.
 If $x=\frac{\boiv}{c_i}$ and $c_i=1$ then again $x\in\ZPhid$.
Finally, if $x=\frac{1}{c_i}\boiv$ and $c_i>1$ then $x$ can't be special.
 Indeed $(\frac{1}{c_i}\boiv,\alpha_i)=1/c_i<1$, thus no member of the family 
of
hyperplanes in $\calh_\Z$ parallel to $H_i$ pass through $x$.
\qed

 \pu
\lemma \label{all special iff an}
{\sl All the vertices of the complex $\Sigma$ are special  if and only if
$\Sigma$ is of type $A_n$.}

\proof
  Since $W_a$ preserves the property of the vertex being special
and since it acts transitively on the set of alcoves,
all the the vertices of $\Sigma$ are special if and only if
all the vertices of the standard alcove
$
A=<0,\frac{1}{c_1}\boov ,\ldots,\frac{1}{c_n}{\bonv}>
$
are special.
 As we have already seen in the proof of the preceding lemma,
the non-special points of this alcove are in one one correspondence
with the numbers $c_1,\ldots,c_n$, that are
strictly greater than 1.
 Thus all the vertexes are special if and only if all the numbers $c_i$ in 
the expression $\tilde{\alpha}=\sum_{1\leq i\leq n }c_i\alpha_i$ are equal 
to 1.
 Now inspecting the tables of the root systems in \cite{bourbaki456}, we conclude that
this  happens only in the case of the root system of type
$A_n$.
\qed

\pu
\label{more subcomplexes}
{\bf More subcomplexes}\par
 Note that an intersection of any family of hyperplanes 
from $\calh_\Z$ or corresponding halfspaces is
a subcomplex of a Euclidean Coxeter complex.
 In particular the line 
$\R~\boiv=\cap_{j\not=i}H_{\alpha_j}$
is a subcomplex.
 Note that for any $m\in\Z,i=1,\ldots,n$ the point $m\boiv/c_i$ is the vertex 
of $\Sigma$.
 Indeed 
$(m\boiv/c_i,\tilde{\alpha})=m$ implies that 
$m\boiv/c_i\in H_{\tilde{\alpha},m}$ and 
$(m\boiv/c_i,\alpha_j)=0~,j\neq i$ implies that
$m\boiv/c_i\in H_{\alpha_j,0}$, hence
$m\boiv/c_i$ is an intersection of $n$  hyperplanes 
$H_{\tilde{\alpha},m},~ H_{\alpha_j,0}~,~j\neq i$.
In particular the line segments $[0,\boiv]\subset\R~\boiv $ are the 
subcomplexes of $\Sigma$. 
 
 Next, the sectors
 $S=S_\Pi := \cap_{\alpha \in \Pi}H_{\alpha}^+ $,
$-S=-S_\Pi := \cap_{\alpha \in \Pi}H_{\alpha}^- $
are subcomplexes as well 
as any of their faces (which are the cones)
$
F=F_\Pi := (\cap_{\alpha \in \Pi'\subseteq \Pi}H_{\alpha})
\cap
(\cap_{\alpha \in \Pi\setminus \Pi'}H_{\alpha}^+).
$
\eject

\section{Ordering Euclidean buildings}
\label{ordering buildings}

\pu
\label{def coxeter ordering}
\definitions
 We will consider {\it special edges} of a  Euclidean Coxeter complex
$\Sigma=\Sigma_\Phi$ of dimension $n$, that is the directed edges 
$e\in\Sigma^{(1)}$ such that the origin $\iota e$ of $e$ is a special
vertex.
 Let $E_s$ be the set of all such edges.
 The typical examples of such  edges are given by  the standard alcove
$
A=<0,\frac{1}{c_1}\boov ,\ldots,\frac{1}{c_n}{\bonv}>
$
 constructed in \S \ref{all special iff an}.
 All the directed edges 
$$
[0,\frac{1}{c_1}\boov] ,\ldots,[0,\frac{1}{c_n}{\bonv}]
$$
are special.
 In some sense any special edge $e$ arrives in this way -- indeed, let 
$\iota e=\alpha\in\ZPhid$, then $e-\alpha$ starts at 0 and there is some simple
system $\Pi$ such that $e-\alpha$ is an edge of the alcove
$A_\Pi$, starting at 0.

   More generally call a directed edge $e$ {\it quasi-special} if it lies on a
line segment $[x,y]$ in $\Sigma^{(1)}$ with special vertices $x,y$.
 An example will be any directed edge lying on the line segment
$[0,\boiv]$ since $0,\boiv$ are special, see \S \ref{special are zphid}.
 This remark implies that  any special edge is quasi-special.
 Note that $\widewa$ leaves  $E_{s}$ invariant as well as the set $E_{qs}$ 
of all quasi-special edges.
(It might be that all the  edges in any Coxeter complex,
and hence in any Euclidean building, are  quasi-special, but the proof of 
this is not in the author's possession.)

 Since  the set of all special vertices on the line $L$ of
$\Sigma^{(1)}$ is discrete in Euclidean topology, we conclude 
that for any quasi-special edge $e$ there is a unique minimal (with 
respect  to inclusion) line segment $[x,y]$ in $\Sigma^{(1)}$ with 
special vertices $x,y$, which contains $e$.

 By an {\it ordering} of $\Sigma$ we mean a function
$\tau:E_{qs}\mapsto \{1,\ldots,n\}, n=\dim \Sigma,$
such that

1)\qua for any alcove $A=<x,x_1,\ldots,x_n>$ with a special vertex
$x$ the
function $\tau$ is bijective on the set of special
edges $\{[x,x_1],\ldots,[x,x_n]\}$,

2)\qua $\tau$ is $\widewa$--equivariant,

3)\qua for any line segment in $[x,y]\subset\Sigma^{(1)}$ with special vertices
$x,y$ the ordering function $\tau$ is constant on a set of directed
quasi-special edges lying on $[x,y]$ and oriented from $x$ to $y$.

 \pu
\remarks
 This resembles the notion of a {\it labelling}
of a Euclidean Coxeter complex, which
means that it is possible to partition the vertices into
$n=\dim\Delta+1$ ``types'', in such a way that each alcove has exactly one 
vertex of each type.
 The labellability of a Euclidean Coxeter complex $\Sigma$ follows from the
fact that the $W_a$--action partitions the vertices into $n$ orbits, and we 
can label $\Sigma$ by associating  one label $i=0,1,\ldots,n$ to each orbit.
 In particular the labelling is $W_a$--invariant.
 There  is one obvious distinction between these two notions --
``ordering orders the directed edges'' and ``labelling labels the vertices''.
 For us it is important that the ordering is invariant under translations 
in the apartments.
 In general there are  translations on  $\Sigma$ which preserve the 
structure of a Coxeter complex but does not belong to $W_a$.

 \pu
\theorem
\label{cox can be ordered}
{\sl Any  Euclidean Coxeter complex~ $\Sigma=\Sigma_\Phi$ can be ordered.
 Moreover an ordering is uniquely defined by an  ordering of a set of all
directed edges of a fixed
alcove starting at some fixed special vertex of alcove.}

\proof
 Consider the set of all pairs $(x,A)$ of {\it based alcoves}
that is alcoves $A$ with a fixed special vertex $x$ of it.

 We wish to prove that the extended Euclidean Weyl group
$\widehat{W_a}=\ZPhid \rtimes W$ acts simply transitively on the set of all 
based alcoves,
that is 
for any pair of based alcoves  $(x,A),(x',A')$ there is exactly one element
$w\in\widehat{W_a}$ which takes $x$ to $x'$ and $A$ to $A'$.
  The set of all special vertices
coincides with the coweight lattice $\ZPhid$, (\S\ref{special are zphid}),
consequently
there is a translation from $\ZPhid$
which takes the special vertex $x$ to the special vertex $x'$, hence we may 
assume that $x=x'$.
 Since $x$ is special $S_W(x)=W$ and the family of hyperplanes $\calh_x$
passing through $x$ define a spherical complex canonically isomorphic
to $\coxsph$.
 The alcoves based at $x$ are in one--one correspondence with sectors
of this spherical complex, thus by transitivity there is $w\in S_W(x)=W$
taking $A$ to $A'$. 

 Now let the element $w\in\widewa=\ZPhid \rtimes W$ fix $(x,A)$.
 The translation $t_x:v\mapsto v+x$ belongs to $\widewa$ by lemma
\ref{special are zphid} and $t_x^{-1}wt_x$ fixes $(0,t_x^{-1}A)$.
 In particular $w'=t_x^{-1}wt_x\in W$ and since $W$ acts simply 
transitively on the set of chambers of $\coxsph$ the element
$w'$ is the identity, hence $w$ is the identity.

 Fix a based alcove $(x,A)$ then we assert that any special edge is 
$\widewa$--conjugate to one and only one such an edge of $(x,A)$ having
$x$ as an origin. 
 Let $e$ be a special edge.
 Since $\ZPhid$ acts simply transitively
on the set $\Sigma^{(0)}_{spec}$ one may assume that $x=0$.
 For  the same reason one may assume that 0 is the origin of $e$.
 Now the sector corresponding to $A$ is a fundamental domain for the 
action of $W$ on the corresponding spherical complex and we can conclude
that $W$ is 
$W$--conjugate to one and only one such  edge of $(x,A)$ having
$x$ as an origin.  
  A priori this does not mean that it  is 
$\widewa$--conjugate to one and only one such edge of $(x,A)$ having
$x$ as an origin.
 But the stabilizer of $0$ in $\widewa$ and that of in $W_a$ is the
same (=$W$).  

 The properties just proven  allow us to order all special edges in
the following way.
 Order the edges of any fixed based alcove $(x,A)$ starting at $x$
arbitrarily and extend the ordering 
in a $\widehat{W_a}$--equivariant way.

 What is left is to extend the ordering to all quasi-special edges.
 Any such  edge $e$ lies on a line segment 
$[x,y]\subset\Sigma^{(1)}$ with special $x,y$ and there is unique minimal
such segment 
(it is fully defined by the condition that there are only two
special vertices on it, namely $x$ and $y$.)
 Let $[x,y]$ be such a segment and let $e$ be oriented from $x$ to $y$.
 If $e$ starts at $x$ then it is a special edge and already has got
its label.
 If not then we assign to $e$ the same order, which has the special edge
on $[x,y]$, beginning at $x$.
 This assigning is a canonical  one thus we have a well defined function
$\tau$ on the set of all quasi-special vertices.
 Let us verify conditions 1)--3) in the definition of an ordering (\S 
\ref{def coxeter ordering}).
 Condition 1) concerns only special edges and 
the bijection desired is given
by the construction.
 $\widewa$--equivariance for special edges again is given by the 
construction and for quasi-special vertices it follows from the 
observation that the minimal $[x,y]$ attached to any such an edge
is obviously $\widewa$--equivariant.
 Finally 3) is given by the construction in the case of minimal $[x,y]$
and clearly any non minimal such a segment is the union of minimal one.
 The only jumping of the order can be in the internal special vertex, but
now we can use the fact that the ordering of 
special edges is $\widewa$--equivariant and in particular $\ZPhid$--equivariant.
\qed

 \pu
\corollary
\label{ordtypcow}
{\sl The ordering of any  Euclidean Coxeter complex corresponding to a root
system is uniquely defined by an
(arbitrary) choice of type function $\tau$ from the set $\{\boiv\}$ of
all fundamental coweights corresponding to some simple system $\Pi$
to the set $\{1,\ldots,n\}.$}

\proof
 Indeed, there is an alcove of the form
$$
A=<0,\frac{1}{c_1}\boov,\ldots,\frac{1}{c_n}\bonv>,
$$
where $c_1,\ldots,c_n$ are defined by the expression of the highest root
$$\tilde{\alpha}=\sum_{1\leq i \leq n}c_i\alpha_i,$$ see \S\ref{nice alcove}.
 The directed edges $[0,\frac{\boiv}{c_i}]$ of this alcove (based in 0)
are in one one correspondence with the fundamental coweights
$\{\boiv\}$.
\qed

\pu
\label{def buildings}
\definitions
  We adopt the following direct definition of a   Euclidean
building, see \cite{brownbuildings}.
 We call a simplicial piecewise Euclidean complex $\Delta$  a
{\it Euclidean building}  if it can expressed as the
union of a family of subcomplexes $\Sigma$, called {\it apartments}
or {\it flats}, satisfying the
following conditions:

    B0)\qua There is a Euclidean Coxeter complex $\Sigma_0$, such that for
each apartment $\Sigma$ there is a simplicial isometry
between $\Sigma$ and $\Sigma_0$.

    B1)\qua Any two simplices of $\Delta$ are contained in an apartment.

    B2)\qua Given two apartments $ \Sigma,\Sigma'$ with a common top-dimensional
simplex (=chamber or alcove), there is an isomorphism
       $\Sigma \stackrel{\approx}{\to} \Sigma'$
        fixing  $\Sigma \cap \Sigma'$  point wise.

Note that isomorphisms in B2) are uniquely defined.
 A Euclidean building has a canonical metric, consistent with the Euclidean
structure on the apartments.
   So each apartment  $E$ of    $\Delta$   is a Euclidean space with a
 metric
\mbox{$\vert x-y\vert_E,x,y\in E$}.
 Moreover, the isomorphisms 
$\Sigma_0 \stackrel{\approx}{\to} \Sigma$ and isomorphisms between  
apartments
given by the building axiom B2) are isometries.
 The metrics
$\vert x-y\vert_E$ can be pieced together to make the entire building
$\Delta$  a metric space. The resulting metric will be denoted by
$x,y \mapsto \vert x-y\vert$. It is known (see 
\cite{brownbuildings}, Theorem VI.3) that the
metric space $\Delta$  is complete.
 Besides for any  $x,y\in \Delta$  the line segment $[x,y]$
is independent of the choice of $E$ and can be characterized by
$$
[x,y] = \lbrace z\in X ; \vert x-y \vert = \vert x-z\vert +
 \vert z-y \vert \rbrace.
$$
 Moreover  $[x,y]$  is {\it geodesic}, that is, it is  the
                shortest path joining  $x$ and $y$
and there is no other geodesic joining $x$ and $y$.

\pu{\bf Local geodesics are geodesic}\par
\label{glg}
 A {\it local geodesic} is  defined as a finite union of line segments
such that the angles between subsequent segments are equal $\pi$.
 The important fact is that in Euclidean buildings  local geodesics
are geodesics.
 This fact is valid for the much more general case
of CAT(0)--spaces, \cite{bridsonhaef}, Proposition II.10.

 \pu
\label{deford}
\definitions
 We extend the definitions from  \S \ref{def coxeter ordering} 
to the case of Euclidean buildings.
 The notion of a special vertex can be defined in the case of a Euclidean
building just by putting the vertex into an apartment.
 This is well defined because any isometrical isomorphism $\phi$ between
apartments takes special vertices to special ones 
(since the same
number of walls pass through $x$ and $\phi(x)$).
 Analogously  the notion of a special edge is well defined.
  
 Call a directed edge $e$ of a building $\Delta$
{\it quasi-special} if it lies on a line segment $[x,y]$ in $\Delta^{(1)}$ 
with special vertices $x,y$.
 Let $E_s, E_{qs}$ be the sets of all special and quasi-special edges 
respectively. 
 Note the inclusion $E_s\subseteq E_{qs}$ which follows from the fact that any
special edge $e$ can be brought to the form 
$e=[0,\frac{1}{c_i}{\boiv}]$ by the extended Weyl group, see
\S \ref{def coxeter ordering}, and 
now it is contained in the line segment $[0,\boiv]$ in $\Delta^{(1)}$ 
which has the special end points.

 By an {\it ordering of an Euclidean building} $\Delta$ 
we mean a function
 $\tau:E_{qs}\mapsto \{1,\ldots,n\}, n=\dim \deltaone,$
such that when restricting to the set of quasi-special edges of any 
apartment it becomes an ordering of a corresponding Coxeter complex.

\pu
\lemma
\label{iso acts on orderings}
{\sl
Let $\Sigma_\Phi,\Sigma'_\Phi$ be two realizations of a Euclidean
Coxeter complex in Euclidean spaces $V,V'$ 
respectively and 
let $\phi:V\rightarrow V'$ be an  isometry taking 
$\Sigma_\Phi$ to $\Sigma'_\Phi.$ 
 If $\tau$ is any ordering on $\Sigma'_\Phi$, then $\tau\phi$ is
an ordering on $\Sigma_\Phi$.
}

\proof
 The only non obvious condition is that $\tau\phi$ is 
$\widewa$--equivariant that is $\tau\phi w=\tau\phi$ for any
$w\in\widewa$.
 But the last equation is equivalent to 
$\tau\phi w\phi^{-1}=\tau$ and the assertion follows from the
fact that  $\phi$ normalizes
the extended Weyl group by the lemma
\ref{weyl groups are invariant}.
\qed

\pu
\theorem
{\sl
 Any  Euclidean building can be ordered.
 Moreover an ordering is uniquely defined by an ordering of special
edges of a fixed alcove of a building starting at some fixed special vertex.
 In particular, there are  only finitely many  orderings on any 
Euclidean  building.
 Isomorphisms in (B2) in \ref{def buildings} can be taken 
to be order-preserving.
 If $\phi$ is an automorphism of a Euclidean building $\Delta$
and $\tau$ is its ordering, then
$\tau\phi$ is  
again an ordering on $\Delta$.}

\proof
 We follow the proof of a labellability of a building \cite{brownbuildings}, 
Chapter  IV, Proposition 1.
 Fix an arbitrary alcove $A=<x,x_1,\ldots,x_n>$ with special vertex
and order the edges $[x,x_1],\ldots,[x,x_n]$ by $1,\ldots,n$ 
respectively.
 If $\Sigma$ is any apartment containing $A$, then as was proved in
\S \ref{cox can be ordered},
there is a unique ordering $\tau_\Sigma$ which agrees
with the chosen ordering on $A$. 
 For any two such apartments $\Sigma,\Sigma'$
the orderings $\tau_\Sigma,\tau_{\Sigma'}$ agree on the
special edges of $\Sigma\cap\Sigma'$;
this follows from the fact that  that $\tau_{\Sigma'}$ can be constructed
as $\tau_{\Sigma} \phi$, where $\phi:\Sigma\rightarrow\Sigma'$  is the 
isomorphism fixing $\Sigma\cap\Sigma'$.
 Since by \S\ref{iso acts on orderings}
$\tau\phi$ is again an ordering and since it coincide on the based
alcove in $\Sigma\cap\Sigma'$ they coincide everywhere.
 The various orderings $\tau_\Sigma$ therefore fit together  to give an 
ordering $\tau$ defined on the union of the apartments containing $A$.
But this union is all of $\Delta$.

 To prove the second assertion,  note that 
isomorphisms $\phi:\Sigma\rightarrow\Sigma'$ fixing $\Sigma\cap\Sigma'$
also fixes some alcove in $\Sigma\cap\Sigma'$ pointwise,
hence it preserves the order.

 Finally, to prove that $\tau\phi$ is again an ordering, just note 
that $\phi$ leaves invariant the set of all  (quasi)-special edges.
 The $\widewa$--invariance follows from the lemma 
\ref{weyl groups are invariant}.
 \qed

 \pu{\bf Standard ordering}\par
\label{standard ordering}
Let us order the fundamental coweights $\boov,\ldots,\bonv$ 
of the corresponding root system as
they are naturally ordered in the tables of root systems given in \cite{bourbaki456}.
 This gives the {\it standard ordering of the standard alcove}, see 
\S\ref{nice alcove}
and thereby the {\it standard ordering of the building}.

\section{Definition of a combing}
\label{sec def combing}
 Let $\Delta$ be an ordered  Euclidean building.
 We wish to  construct  a combing $\calc$ on $\Delta$ which 
consists of edge paths in the 1--skeleton $\Delta^{(1)}$ and is 
$\widewa$--equivariant when restricted to any apartment.
 By definition an {\it edge path} $\alpha$ in a graph $\Delta^{(1)}$
is a map of  the interval $[0,N]\subset\N$
into $\Delta^{(0)}$
such that $\forall i\in[0,N-1]$ the vertices  $\alpha (i),\alpha(i+1)$ are the 
end points of an edge in $\Delta^{(1)}$.
It is convenient to consider $\alpha$ as an ultimately constant map
by extending it to a map of $[0,\infty ]$ making
$\alpha(t)$ stop after $t=N$, ie, by setting 
$\alpha(t)=\alpha(N)$ for $t\geq N$.

\pu
\label{def combing}
{\bf Combing the building}\par
 Take any  special vertices $x,y\in \deltaone$ and put them into some 
apartment $\Sigma$.
 We consider $\Sigma$ as a vector space, taking $x$ as an origin.
 Since $x$ is a special vertex, the walls $H\in\calh_\Z$,
passing through $x$, define the structure of a spherical Coxeter complex
$\Sigma_{sph}(x)$.
 In particular 
$y$ lies in some closed sector $S$ of this spherical complex, based in $x$.
 To  $S$ one can canonically associate a based alcove $(x,A_S)$
-- this is a unique alcove in $S$ with $x$ as one of its vertices.
 As with any sector, $S$ is a simplicial cone and the rays defining it are
the rays spanned by the edges of the alcove $A_S$, started in
$x$.
 More exactly, let $A_S=<\xon>$  with the edges $[x,x_i]=e_i$ of  type 
$i, i=1,\ldots,n$.
 These edges constitute the basis of a vector space structure on $\Sigma$,
corresponding to a choice $x$ as an origin.
 For any $e_i$ we define $e_i$--direction as the set of all rays of the
form $z+\R e_i, z\in\Sigma$.
 We have $S=\sum_{1\leq j\leq n}\R_+ e_j$ relative to our vector space 
structure.
 Now $y\in S$, hence $y=\sum_{1\leq i\leq n} m_ie_i$ with  
$m_1,\ldots,m_n\geq 0.$
 We are able now to define the path $\alpha_{xy}$ of a combing
$\calc$ connecting $x$
(= 0 relative to a vector space structure) and $y$.
 This is a concatenation of the line segments
\begin{equation}
[0,m_1e_1],[m_1e_1,m_1e_1+m_2e_2], \ldots,
[\sum_{1\leq i\leq n-1} m_ie_i,\sum_{1\leq i\leq n} m_ie_i],
\end{equation}
passing in this order.
  Geometrically speaking, $\alpha_{xy}$ is a concatenation of an ordered 
sequence
of $n$ segments (degenerate segments are allowed), such that the $i$-th 
segment is parallel to the line $\R~ e_i$ (degenerate segment is considered as
to be parallel to any line).

  Define now  a combing $\calc$ by collecting all the paths of the form
$\alpha_{xy}$ (for the special vertices $x,y$) as well as all
their  prefix subpaths in $\deltaone$.

 \pu
\label{special graph}
{\bf Graph structure $\deltaonespec$}\par
 Apart from the natural simplicial graph structure on $\Delta^{(1)}$ we 
wish to use another rougher simplicial graph structure $\deltaonespec$.
 The vertices of $\Delta^{(1)}_{spec}$ are the special vertices
of $\Delta^{(1)}$.
 Two special vertices $x,y$ are connected by the edge $[x,y]$
(which is a line segment joining these vertices)
if $[x,y]$ lies in $\deltaone$ and there are no other special vertices
between $x$ and $y$.
 The main example of an edge in $\Delta^{(1)}_{spec}$ will be 
the line segment $[x,x+\boiv], i=1,\ldots,n$ in a Euclidean Coxeter complex, 
where $x$ is a special vertex,
see \S\ref{special are zphid}.
 One obtains other examples  from the observation that the extended
Weyl group acts preserving the structure $\deltaonespec$.

 Let $\calcspec$ be a subcombing of $\calc$ consisting of those paths 
that connect only the special vertices of $\Delta$.
 Note that a combing $\calc_{spec}$ 
gives rise naturally to a combing of a graph $\deltaonespec$
which we will denote by the same symbol.
\pu
\label{qi1}
\lemma
{\sl
The natural embeddings $\deltaonespec\subseteq\deltaone\subset\Delta$ induce
quasi-isometry between the graphs  $\deltaonespec$, $\deltaone$
with their graph metrics and the building $\Delta$ with 
its piecewise Euclidean metric.}

\proof
 Recall some definitions (see  \cite{neumannshapiroboundaries}).
Given $\lambda\geq 1$ and $\epsilon\geq 0$, a map 
$f:X\rightarrow Y$ of
metric spaces is a $(\lambda,\epsilon)$--{\it quasi-isometric map} if
$$
\oneoverlambda d_Y(x,y)-\epsilon\leq d_X(f(x),f(y))\leq 
\lambda d_Y(x,y)+\epsilon 
$$
for all $x,y\in X$. 
 If $X$ is an interval, we speak of a
{\it quasigeodesic path} in $Y$. 
 Two metric spaces $X$ and $Y$ are
{\it quasi-isometric} if there exists a quasi-isometric map 
$f:X\rightarrow Y$ such that $Y$ is a bounded neighborhood of the image of 
$f$.
 Then $f$ is called a {\it quasi-isometry}.
 If the constant $\epsilon$ above is zero, then we speak about
{\it Lipschitz maps} and {\it Lipschitz equivalence}.

 Firstly  we prove that the embedding  
$\deltaonespec\subset\Delta$ is a Lipschitz equivalence.
 If $\lambda$ is the Euclidean length of the longest  edge
in $\deltaonespec$, then clearly $d_\Delta\leq\lambda d_{1,sp}$, where 
$d_{1,sp}$ is a graph metric on $\deltaonespec$.
 On the other hand, if $x,y$ are the special vertices, then by definition 
$\alpha_{xy}$ is a concatenation of the line segments
$$
[0,k_1\boov],[k_1 \boov,
k_1 \boov+k_2 \botv], \ldots,
[\sum_{1\leq i\leq n-1} k_i\boiv,\sum_{1\leq i\leq n}
k_i\boiv]
$$
passing in this order.
 Its graph length in $\deltaonespec$ is equal to 
$\sum_{1\leq i\leq n}k_i$ and this is an $\ell^1$--distance between
$x,y$ in $\ell^1$--metric on
$\Sigma$, given in coordinates attached to a basis $\boov,\ldots,\bonv$.
 Since, up to translations in $\Sigma$, there are only finite number
of bases of this type and since the $\ell^1$--metric is Lipschitz equivalent
to the Euclidean metric there is a constant $\lambda_2>0$, such that
 $d_{1,sp}(x,y)\leq \lambda_2 d_\Delta(x,y)$ -- this proves that
an embedding $\deltaonespec\subset\Delta$ is a Lipschitz equivalence.

 Now, any edge in $\deltaonespec$ consists of several edges of
$\deltaone$ and let $\lambda_3$ be the largest number of the edges in 
$\deltaone$ lying on the edge of $\deltaonespec$.
 Then for the graph metric $d_1$ on $\deltaone$
we have an inequality $d_1\leq \lambda_3 d_{1,sp}$.
 Conversely, if $\lambda_4$ is the Euclidean length of the longest edge
of the graph $\deltaone$, then $d_\Delta\leq\lambda_4 d_1$ and since
$d_{1,sp}$ and $d_\Delta$ are Lipschitz equivalent we get the 
inequality $d_{1,spec}\leq\lambda_5 d_1$ for a suitable 
positive constant $\lambda_5$.
\qed

\pu
\label{qi2}
\lemma
{\sl 
The induced metric on 
$\deltaonespec$ as well as on $\deltaone$
(induced from the Euclidean metric on a building)
is Lipschitz equivalent to the edge path metric on these  graph.}

\proof
 This simple observation is indeed true in a more general
situation of $\R$--{\it graphs}, that is simplicial graphs,
in which any edge is endowed with  a metric, making it to be isometric
to a segment of a real line.
 Namely if the lengths of the edges in such a graph are bounded from 
above and from
below by some positive constants, then the $\R$--metric on a graph is 
Lipschitz equivalent to the edge path metric on this  graph. 
 This is clear.
 Now in our situation of Euclidean Coxeter complexes there
are only finitely many isometry types of the edges, so we have
the bounds just mentioned.
\qed

\ppar
{\bf Properties of $\calc$} 

\pu
\label{combing is quasigeodesic}
{\sl  $\alpha_{xy}$ is quasigeodesic relative to the
edge path metric on $\Delta^{(1)}.$}\par
 Obviously the path $\alpha_{xy}$ is geodesic relative to an 
$\ell^1$--metric on
$\Sigma$, given in coordinates attached to a basis $e_1,\ldots,e_n$.
 This $\ell^1$--metric is Lipschitz equivalent to a Euclidean metric on 
$\Sigma$, which is Lipschitz equivalent to a graph metric on $\Sigma$.
 These equivalences preserve the the quasigeodesicity of paths, 
hence our path is quasigeodesic.
 (It seems likely that indeed $\alpha_{xy}$ is geodesic).

  \pu
{\sl Support of the path}\par
 The path $\alpha_{xy}$ uniquely defines a cone 
$C_\alpha=\sum\{\R_+e_i;m_i>0\}$ which  is called a {\it support} of
$\alpha$. 
(It  is uniquely defined in $\Sigma$, not in $\Delta$!). 
 $C_\alpha$ is a subcomplex of $\Sigma$, see \S\ref{more subcomplexes}.
 Note that $\alpha$ travels inside of $C_\alpha$ (indeed it is
the least cone in which it travels, since it can be
defined as the cone  span of $\alpha$).
 Note also that $C_\alpha$ is the smallest face of $S$,
containing $y$ in its interior.
 In particular, any sector, containing $y$, contains also $C_\alpha$.

 \pu\label{edge path}
{\sl $\alpha_{xy}$  is an edge path in the 1--skeleton $\Delta^{(1)}$.}\par
 Let $S=S_\Pi$ for some simple system $\Pi=\{\alphai\}$ of an 
underlying root system.
 Then the corresponding alcove $A_\Pi$ is a closed simplex with the
vertices $0$ and $\boiv/c_i, i=1,\ldots,n$, see \S{\ref{nice alcove}}.
 Here  $\{\boiv\}$ is the dual basis for $\{\alpha_i\}$ (the basis of
coweight lattice $\ZPhid$) and
$\tilde{\alpha}=\sum_{1\leq i \leq n}c_i\alpha_i$ is the corresponding highest
root.
 Hence, after appropriate re-ordering of $\alpha_i$, we may assume that
$e_i=\boiv/c_i, i=1,\ldots,n$.
 The fundamental coweight $\boiv$ is a special vertex 
and lies on the line passing through $e_i$, thus
the line segment $[0,\boiv]$ is an edge path in $\Sigma^{(1)}$.
 Since $y$ is a special vertex it belongs to
the coweight lattice $\ZPhid$ by \S{\ref{special}}, and thus, in the 
notations of \S\ref{def combing}, all the numbers
$k_i=\frac{m_i}{c_i}$ are integers.
 By definition $\alpha_{xy}$ is a concatenation of the line segments
$$
[0,k_1\boov],[k_1 \boov,
k_1 \boov+k_2 \botv], \ldots,
[\sum_{1\leq i\leq n-1} k_i\boiv,\sum_{1\leq i\leq n}
k_i\boiv]
$$
passing in this order.
 We conclude from this formula that the path is a concatenation of the paths,
obtained from the line segments  $[0,\boiv]$ (which are edge paths) by the 
action of $\ZPhid$ and the last
action preserves the simplicial structure on $\Sigma$.

 \pu\label{path in a convex hull}
{\sl $\alpha_{xy}$ is contained in the convex hull of the set $\{x,y\}$.}\par
 By a {\it convex hull} of an arbitrary set  in the building $\Delta$ we 
mean the smallest convex subcomplex containing this set.
 For example apartments are convex, hence the convex hull of any set is
contained in any apartment, which contains this set.
 In a case of a set consisting of two special vertices $x,y$ one can describe 
the convex hull $\chxy$ as a parallelepiped, spanned by $\{x,y\}$.
 More exactly, in the notations of \S\ref{edge path}, we assert that
\begin{equation}
\label{eqparallelepiped}
\chxy =\{z=\sum z_ie_i:~0\leq z_i\leq m_i,~i=1,\ldots,n\}.
\end{equation}
It immediately follows from this description  that $\alpha_{xy}$ is 
contained in  $\chxy$.
 To prove the equality (2), note firstly that
the sector $ S$ is a convex subcomplex as well as $- S$ 
(\S\ref{more subcomplexes}).
 We see immediately that the parallelepiped $P$ on the right hand side is the 
intersection of the sectors $ S$ and $-S+y$ and thus is a convex subcomplex.
 Hence $P\supseteq\chxy$.
 Suppose that $P\not=\chxy$, then $\chxy$ is a proper convex subcomplex of $P$,
containing $x,y$.
 Because $\chxy$ is an intersection of closed half-spaces bounded by elements
of $\calh_\Z$, there is  a half-space $H^+$ defined by some
$H\in\calh_\Z$, containing $\chxy$,
but not $P$.
 (One can prove this fact by showing firstly that $\chxy$ is a convex hull of
finite number of vertices and then follow the standard proof that the convex
polygon is an intersection of half-subspaces, supported on codimension one
faces.)
 Suppose, for instance, that $y$ is not farther from $H$ than $x$.
 Take the hyperplane $H_1\in\calh_\Z$ parallel to $H$ and passing through $y$ 
and let $H_1^+$ be a half-space bounded by $H_1$ and contained in $H^+$.
 Since $H_1$ pass through the special vertex $y$ it belongs to the 
structure $\Sigma_{sph}(x)$ of  spherical
complex in $y$, see \S\ref{def combing}, and since $H_1^+$ contains $x$ it 
contains also the support $C_{\alpha^{-1}}$ of the path $\alpha^{-1}$, 
inverse to $\alpha$.
 But  $C_{\alpha^{-1}}=-C_\alpha +y$, hence it contains the parallelepiped
$P$, contradiction.
 
 This proof does not work in the case when $y$ is not special, since
the set of hyperplanes from $\calh_\Z$ passing through $y$ does 
not constitute the structure of spherical complex.
 But still we can prove that the parallelepiped $P$ in 
(\ref{eqparallelepiped}) is contained in $\chxy$.
 Suppose the contrary, that $P$ is not contained in $\chxy$, then again there 
is  a half-space $H^+$ defined by some
$H\in\calh_\Z$, containing $\chxy$, but not $P$.
 Suppose that 
$y$ is not farther from $H$ than $x$
 (the opposite case was already treated above).
 Define the structure $\Sigma_{sph}(y)$ by translating
such a structure  from any special vertex.
 Relative to this structure $x$ lies in the support 
$C_{\alpha^{-1}}=-C_\alpha +y$, where $\alpha^{-1}$ is the path
inverse to $\alpha$.
 Consider  the hyperplane $H_1\in\calh_\Z$ parallel to $H$ and passing 
through $y$ and let $H_1^+$ be a half-space bounded by $H_1$ and contained in 
$H^+$.
 Since $H_1^+$ contains $x$ it 
contains also  $C_{\alpha^{-1}}=-C_\alpha +y$, hence it contains the
parallelepiped $P$, contradiction.

 \pu
{\sl $\calc$ is $\widewa$--invariant.}\par
 This immediately follows from the fact that the ordering is 
$\widewa$--invariant
and from the geometric interpretation of the paths in $\calc$ just given.

 \pu
{\sl If $x,y$ are  special vertices, then $\alpha_{xy}$  is uniquely 
defined by $x,y$.}\par
 Firstly, the convex hull $\chxy$ is uniquely defined.
 It is a parallelepiped and its 1--dimensional faces 
are ordered.
 Now $\alpha_{xy}$ is a  unique edge path from $x$ to $y$ in $\chxy$
which is a concatenation of 1--faces of $\chxy$ passing in the increasing order.

\section{Fellow traveller property}\label{ftp}

\pu
\theorem
{\sl The combing $\calc$ of an ordered Euclidean building 
$\Delta$ constructed in \S\ref{sec def combing} satisfies
the ``fellow traveller property'',
namely there is  $k>0$ such  that  if
$\alpha,\beta\in\calc$
begin and end   at a distance at most one apart, then
$$
d_1(\alpha(t),\beta (t))\leq k
$$
for all \ $t\geq 0$.
 (The metric $d_1$ is the graph metric on $\deltaone$.)
 The same is true for the combing $\calcspec$ on a graph 
$\deltaonespec$, see \ref{special graph}.}

\proof
A)\qua
 Let us firstly consider the case where  $\alpha,\beta$ begin 
at the same vertex
and end  at a distance  one apart, that is
$[\alpha(\infty),\beta(\infty)]$ is an edge.
 It is easily seen from a $k$--fellow traveller property that for any $c>0$ 
if $\gamma,\gamma'\in\calc$
begin and end at a distance at most $c$  apart, then
$
\vert\gamma(t)-\gamma'(t)\vert\leq kc
$
for all \ $t\geq 0$.
 Note also that we can work in one apartment since the initial vertex of the
paths and the edge $[\alpha(\infty),\beta(\infty)]$ are contained in some
apartment and thus the whole paths lie in this apartment, see 
\S\ref{path in a convex hull}.
 Thus we may assume that $\alpha,\beta$ are contained in the Coxeter complex
$\Sigma$ and start at $0$.
 Associated to $\alpha,\beta$ are their supports $C_\alpha,C_\beta$
in which they travel respectively.
 The intersection $K=C_\alpha\cap C_\beta$ is a simplicial cone
of the form $K= \sum\R_+v_j$ for some set $\{v_j\}$ of special edges.
 Hence  
$$
C_\alpha=\sum\R_+u_i + \sum\R_+v_j
$$
and
$$
C_\beta=\sum\R_+w_k + \sum\R_+v_j
$$
where $\{u_i\},\{v_j\},\{w_k\}$
are the sets of special edges (possibly empty)
and the sets $\{u_i\},\{w_k\}$ do not intersect.
 We will argue by induction on the sum
$
\dim C_\alpha +\dim C_\beta.
$
 The least nontrivial case is when the sum is equal to 2 and both of
$C_\alpha,C_\beta$ are one dimensional, that is they are simplicial
rays.
 If $\alpha$ and $\beta$  have the same direction then obviously they
1--fellow travel each other.
 If not then they diverge linearly with  a speed bounded from below
by a constant not depending on the paths
 (indeed, there only finitely many of possibilities for the angle
between  $C_\alpha$ and $C_\beta$).  
 Thus they could end at a distance one  apart only when they passed
a bounded distance, thus they $k$--fellow travel for some $k>0$.
 A similar argument applies when  $C_\alpha\cap C_\beta=0.$
 The main case is when
$K=C_\alpha\cap C_\beta=\sum\R_+v_j$ is nonzero.
 Again in this case  the argument similar to the above
shows that $\alpha$ (resp.\ $\beta$) 
can move in the $u$--direction (resp.\ in the $w$--direction)
only for a bounded amount of time, $c$ say. 
 The rest of the proof is the reduction to the  case when the paths 
lie in $K$ -- then we can apply inductive hypothesis.
 But we need a definition.
 Consider a path $\gamma,$ consising of two linear subpaths $\gamma',\gamma''$ 
passing in this order and having the directions $v,w$ correspondingly.
 Construct the path $\tilde\gamma$, which goes distance $|\gamma''|$ in 
the $w$--direction
and distance $|\gamma'|$ in the $v$--direction.
 We call $\tilde\gamma$ an {\it  elementary transformation of the 2--portion 
path} $\gamma$.
 Note now that if $|\gamma''|\leq c$, then $\tilde\gamma$ and $\gamma$
~~$2c$--travel
each other.
 Using the elementary transformations as above we can push out all the 
linear
subpaths  of $\alpha$ with $u$--directions to the end of the path.
 So we assume that $\alpha$ and $\beta$ are such from the beginning.
 Of course this operation takes $\alpha$ from the combing $\calc$, but
note that the initial $v$--portion of $\alpha$ continues to lie in $\calc$,
since the elementary transformations do not change the ordering of
$v$--directions.
 Cutting out the $u$--tails and $w$--tails of $\alpha$ and $\beta$
correspondingly, we may assume that $\alpha$ and $\beta$ lie in the
cone $K$ and end within at most distance $c$ for some universal constant
$c$.
 Now let $v_{j_0}$ be of the smallest order in the set $\{v_j\}$.
 Then each of the paths $\alpha,\beta$ move some nonzero time in
$v_{j_0}$--direction, hence they coincide during some nonzero
time.
 Cutting out the longest coinciding part of $\alpha$ and $\beta$ we
may assume that one of them does not contain the 
$v_{j_0}$--direction at all.
 But now the dimension of the support either of  $\alpha$ or $\beta$ 
decreases and we may apply an induction hypothesis.

B)\qua Consider now the general case where
$\alpha,\beta\in\calc$
begin and end  at a distance at most one apart.
 Adding a bounded number of edges to $\alpha$ and $\beta$ one may assume
that  $\alpha$ and $\beta$ end at the special vertices which are
distance at most $c$ apart for some constant $c$, depending only on $\Delta.$
  Drawing a path from $\calc$ connecting $\alpha(0) $ and $\beta(\infty)$
and making use the part A) of the proof we reduce the problem to  the
case when $\alpha,\beta$ end at the same special vertex.
 Denote by $\calc_{spec}$ the set of all paths from $\calc$ connecting
special vertices. 
 Now consider the combing $-\calc_{spec}$ consisting by  definition
of the paths which are inverse to the paths from $\calc_{spec}$.
 Obviously $-\calc_{spec}$  is obtained in the same way as 
$\calc_{spec}$ but reversing the underlying ordering.
 Thus  for some $k'$  the combing $-\calc_{spec}$ satisfies 
the $k'$--fellow 
traveller property for paths which begin at the same vertex.
 Now let $\alpha^{-1},\beta^{-1}$ be the paths inverse
to  $\alpha,\beta$ so that $\alpha^{-1},\beta^{-1}\in-\calc_{spec}$.
 Consequently $\alpha^{-1},\beta^{-1}$ $k'$--fellow travel each other.
This doesn't imply that $\alpha$ and $\beta$ are $k'$--fellow
travellers since they arrive at $\alpha(0)$ and $\beta(0)$ at different
times.
 Let $N_\alpha$ and $N_\beta$ the the length of $\alpha$
and $\beta$ respectively.
 Then $\alpha^{-1}(t)=\alpha(N_\alpha -t)$ (assuming that $\alpha$
extended to the negative times in an obvious way).
 We have $|\alpha(t)-\beta(t)|=|\alpha^{-1}(N_\alpha-t)-\beta^{-1}(N_\beta-t)|
\leq 
|\alpha^{-1}(N_\alpha-t)-\beta^{-1}(N_\alpha-t)|+ 
|\beta^{-1}(N_\alpha-t)-\beta^{-1}(N_\beta-t)|.$
 The first modulus is bounded since $\alpha^{-1},\beta^{-1}$
fellow travel each other.
 The second modulus is bounded since the paths
$\beta^{-1}(N_\alpha-t), \beta^{-1}(N_\beta-t)|$
differ  by only a bounded time shift.

C)\qua  Let 
$\alpha,\beta\in\calc$ be the paths from a combing $\calcspec$, 
beginning  and ending   at a distance at most one apart.
 There is a constant $c>0$, depending only on $\Delta$, such that
any edge in $\deltaonespec$ is of a length $\leq c$.
 Then, applying part B) we prove that $\alpha,\beta$ fellow travel each
other relative to $\deltaone$--metric and since the metrics 
$d_1, d_{1,spec}$ are Lipschitz equivalent, we get the fellow traveller 
property for the combing $\calcspec$  on a graph $\deltaonespec$.
\qed

\section{Recursiveness of a combing $\calc$}
\label{recursive}
 In this section $\Delta$ will be an ordered Euclidean building
with a standard ordering, see \S\ref{standard ordering}.

 The definition of  $\calc$ given above is "global" in the sense that
a path from  $\calc$ "knows" where it goes to.
 In this section we show that a path from $\calc$ can be defined 
by a simple local ``direction set'':
namely any  pair of consecutive directed  quasi-special
edges $\{e_1,e_2\}$ shall be one of the following two types:

1)\qua $\{e_1,e_2\}$ is {\it straight}, that is the angle 
between $e_1$ and $e_2$ is equal $\pi$
and hence the union $e_1\cup e_2$ is the line
segment of length 2 in the edge path metric,

or

2)\qua the type $i$ of $e_1$ is strictly less than the type $j$ of $e_2$,
the end of $e_1$
(= the origin of $e_2$) is a special vertex 
and $(e_1,e_2)=(\frac{1}{c_i}\boiv,\frac{1}{c_j}\bojv).$

 Define  $\calc'$ to be the family of all paths $\gamma$ in $\deltaone$, 
in which  any pair of consecutive edges satisfy  
either 1) or 2).

 \pu
\theorem
{\sl If $\Delta$ is a Euclidean building  of one of the following three types
$A_n,B_n,C_n$, then the  combings  $\calc$ and $\calc'$ coincide.}

\proof
It follows immediately from the global definition of $\calc$
that  it is contained in  $\calc'$.
 The proof of the converse proceeds  by induction on the number of line 
segments constituting the path $\gamma\in\calc'$.
 Take $\gamma \in \calc'$ and write it as
$\gamma=\beta\cup e$, where $e$ is the last edge of $\gamma$ and $\beta$ is
the portion of $\alpha$, preceding $e$ .
 By the induction hypothesis $\beta\in\calc$ thus
$\beta\subset\chxy$, where $x=\beta(0),y=\beta(\infty))$ in view of 
 \S {\ref{path in a convex hull}}.
 Let $\Sigma$ be an apartment containing both $\beta(0)$ and $e$.
 Since $\Sigma$ is convex it contains $\chxy$
and thereby $\beta$.
 Consequently  all our path $\gamma$ is contained in $\Sigma$.
 Take $x$ as an origin and identify $\Sigma$ with the standard
Coxeter complex $\Sigma_\Phi$.
 Thus we may assume that $\beta $
lies in a standard sector $\sum\R_+\boiv$ and
$$
y=
m_1\frac{\boiov}{c_{i_1}}+\ldots+
m_r\frac{\boirv}{c_{i_r}},
$$
where all the coefficients $m_i$ are natural numbers.
 Let $e_1$ be the last edge  of $\beta$, then it  is parallel to $\boirv$ and 
is of  type $i_r$.
 With the notation $j=i_r$ and
 by definition of $\calc'$ the type $k$ of $e$ is not smaller than $j$
and the pair
 $\{e_1,e\}$
is one of the following two types:

1)\qua $\{e_1,e\}$ is {\it straight}, that is the angle between these two
vectors is zero and, hence, by \ref{glg} the 
union $e_1\cup e$ is the line
segment of the length 2 in edge path metric

or

2)\qua the type $j$ of $e_1$ is strictly less than the type $k$ of $e$
and
$(e_1,e)=$\break$(\frac{1}{c_i}\boiv,\frac{1}{c_j}\bojv).$

 Since $y$ is a special vertex its stabilizer 
$W(y)$ is conjugate to $W$ -- the  Weyl group of a root system $\Phi$.
 Since the set of all the edges of type $k$ starting in 
$\beta(\infty)=y$ is an orbit $W\bokv$
we have $e=w\bokvnack$ for some $w\in W(y).$
 If one could find this $w\in W$ in such a way that it fixes all 
$\boiv,i\leq j,$ then such $w$ fixes $\beta$,
 since $\beta$
lies in a Euclidean subspace spanned by the vectors $\boiv,i\leq j$.
 Now applying $w^{-1}$ to the path $\alpha$, we get that
 $w^{-1}\alpha=\beta\cup \frac{1}{c_k}\bokv$, hence  
$w^{-1}\alpha\in\calc$.
 Taking into account
that $W_a$ preserves $\calc$, we get that $\alpha\in\calc.$
 The problem now is to find $w\in W(y)$
with the properties as above.
 As was mentioned above the last edge of $\beta$ is parallel to $\bojv$ and 
is of  type $j$.
 Again by definition of $\calc'$ we deduce that
$
(\onenacj{\bojv}, e)=
(\onenacj\bojv,\onenack\bokv).
$

 Thus, to finish the proof we need the following technical lemma.

 \pu
\label{technical lemma}
\lemma
{\sl
{\rm(a)}\qua
 Let $\Phi$ be a root system of one of the  types $A_n,B_n,C_n$, 
given by the tables in \cite{bourbaki456}, pages 250--275.
 Order fundamental coweights by their indices as they are given in 
\cite{bourbaki456}.
 Let $\omega$ be a coweight of type $k$ and
$$
(\bojv, \omega)=(\bojv,\bokv)
$$
for some  $j< k$.
 Then there is $w\in W$ fixing  all the vectors $\boiv,i\leq j$ and such
that $\omega=w\bokv$. 

{\rm(b)}\qua The assertion is not true for the remaining classical
case when $\Phi$ is of type $D_n$.}

 \pu
{\bf Proof of lemma \ref{technical lemma}}

{\bf  (a)\qua Case $A_n, n\geq 1$}

 Denote by $\epsilon_0,\ldots,\epsilon_n$ the standard basis of
$\Rnone,n\geq 1$.
 Let $V$ be the hyperplane  in $\Rnone$ consisting of vectors whose
coordinates add up to $0$.
 Define $\Phi$ to be the set of all vectors  of squared length 2 in the
intersection of $V$ with the standard lattice
$\Z\epsilon_0 +\ldots +\Z\epsilon_n$.
 Then $\Phi$ consists of the $n(n+1)$ vectors:
$$
\epsi -\epsj , ~~\ 0\leq i\not= j\leq n
$$
and $W$ acts as a permutation group $S_{n+1}$ on  basis
$\epsilon_0,\ldots,\epsilon_n$.

 For the simple system $\Pi$ take
$$
\alpha_1=\epsilon_0-\epsilon_1,
\alpha_2=\epsilon_1-\epsilon_2,\ldots,
\alpha_n=\eps_{n-1}-\epsn .
$$
Then the highest root is 
$$
\tilde{\alpha}=\epszero-\epsn =
\alpha_1+\alpha_2+\ldots+\alpha_n .
$$
The fundamental coweights are
$$
\bojv=(\epszero+\ldots+\eps_{j-1})-\frac{j}{n+1}
\sum_{i=0}^{n}\epsi ,
~~
1\leq j\leq n~ .
$$
Let the coweight $\omega$ satisfies the hypotheses of the lemma.
 Then since the type is $W$--invariant  $\omega=u\bokv$ for some $u\in W$ .
 As  $W$ acts by permutations on the basis $\boiv$
\begin{equation}\label{e}
\omega=\sum_{1\leq r\leq k}\epsir -\frac{k}{n+1}
\sum_{i=0}^{n}\epsi .
\end{equation}
By the hypotheses of the lemma
\begin{equation}\label{oeeqoo}
(\bojv, \omega)=(\bojv,\bokv)
\end{equation}
for some  $j<k$.

The left hand side of (\ref{oeeqoo}) is
$$
\Big(
\epszero+\ldots+\eps_{j-1}-\frac{j}{n+1}
\sum_{i=0}^{n}\epsi,~
\sum_{1\leq r\leq k}\epsir -\frac{k}{n+1}
\sum_{i=0}^{n}\epsi
\Big)
=
$$
\begin{equation}\label{anleft}
\mbox{card}\{r|i_r\leq j-1\}-\frac{jk}{n+1}.
\end{equation}
The right hand side of (\ref{oeeqoo}) is
$$
(\bojv,\bokv)=
\Big(\epszero+\ldots+\eps_{j-1}-\frac{j}{n+1}\sum_{i=0}^{n}\epsi ,~
\epszero+\ldots+\eps_{k-1}-\frac{k}{n+1}
\sum_{i=0}^{n}\epsi
\Big)=
$$
\begin{equation}\label{anright}
\frac{j(n+1-k)}{n+1}=j-\frac{jk}{n+1}.
\end{equation}
Now comparing (\ref{anleft}) and (\ref{anright}), we conclude that
$$
\mbox{card}\{r|i_r\leq j-1\}=j,
$$
hence from (\ref{e})
$$
e=(\epszero+\ldots+\eps_{j-1})+
\epsij+\ldots+\epsik
-\frac{k}{n+1}
\sum_{i=0}^{n}\epsi.
$$
 But $\omega=u~\bokv$ and
$$
\bokv=(\epszero+\ldots+\eps_k)-\frac{k}{n+1}
\sum_{i=0}^{n}\epsi,
$$
consequently $\omega=w\bokv$ for some $w\in W$ fixing
all the vectors $\eps_0,\epsone,\ldots,\eps_{j-1}.$
 Since $\boiv, i\leq j$
are the linear combinations of the vectors
$$\eps_0,\epsone,\ldots,\eps_{j-1},
\frac{j}{n+1}\sum_{i=0}^{n}\epsi$$ (the last one is fixed by $W$), we get
that  $\boiv, i\leq j$ are also fixed by $w$.

\bigskip

{\bf Case $B_n, n\geq 2$}

 Denote by $\epsilon_1,\ldots,\epsilon_n$ the standard basis of $\Rn,n\geq 2$.
 Define $\Phi$ to be the set of all vectors  of squared length 1 or 2 in the
standard lattice
$\Z\epsilon_1 +\ldots +\Z\epsilon_n$.
 So $\Phi$ consists of the $2n$ short roots $\pm\epsi$
and the  $2n(n-1)$ long roots
$\pm\epsi \pm\epsj , ~~(i< j)$, totalling $2n^2.$
 For the simple system $\Pi$ take
$$
\alpha_1=\epsilon_1-\epsilon_2,
\alpha_2=\epsilon_2-\epsilon_3,\ldots,
\alpha_{n-1}=\eps_{n-1}-\epsn,
\alpha_n=\epsn .
$$
Then the highest root
$$
\tilde{\alpha}=\epsilon_1+\epsilon_2 .
$$
 The Weyl group $W$ is the semidirect product of $S_n$ (which permutes
$\epsi$) and $(\Z/2)^n$ (acting by sign
changes on the $\epsi$), the latter
normal in $W$.

 The fundamental coweights are
$$
\boiv=\epsilon_1+\epsilon_2+\ldots+\epsi ,
~~
1\leq i\leq n~ .
$$
Let the coweight $e$ satisfies the hypotheses of the lemma.
 Then since $W\bokv$ is the set of all edges of type  $k$
we have $\omega=u\bokv$ for some $u\in W$ .
 As $W$ acts by ``sign'' permutations on the basis $\boi$ we have that
\begin{equation}\label{be}
\omega=\sum_{1\leq r\leq k}\pm\epsir .
\end{equation}
By hypotheses of the lemma
\begin{equation}\label{boeeqoo}
\bojv \omega=\bojv\bokv
\end{equation}
for some  $j<k$.

The left hand side of (\ref{boeeqoo})  is equal to
\begin{equation}\label{bleft}
(
\epsilon_1+\ldots+\epsj
)
(
\sum_{1\leq r\leq k}\pm\epsir
)
\leq
\mbox{card}\{r|i_r\leq j\}.
\end{equation}
The right hand side of (\ref{boeeqoo}) is equal
\begin{equation}\label{bright}
\bojv\bokv=j.
\end{equation}
Now comparing (\ref{bleft}) and (\ref{bright}), we conclude that
$$
\mbox{card}\{r|i_r\leq j\}=j,
$$
hence from (\ref{be})
$$
\omega=(\pm\epszero\pm\ldots\pm\eps_j)+
\pm\epsij\pm\ldots\pm\epsik.
$$
 But $\omega=u~\bokv$ and
$$
\bokv=\epsilon_1+\epsilon_2+\ldots+\epsk ,
~~
1\leq k\leq n~ .
$$
Since $\boiv, i\leq j$
are  linear combinations of the vectors
$\eps_0,\epsone,\ldots,\eps_{j-1},$  we get
that  $\boiv, i\leq j$ are also fixed by $w$.

 But $\omega=w~\bojv$, consequently $w$ can be chosen in such a way that
it fixes the vectors $\epsone,\ldots,\eps_j.$

\bigskip
{\bf  Case $C_n, n\geq 2$}

 Starting with $B_n$, one can define $C_n$ to be the
inverse root system.
 It consists of the $2n$ long roots $\pm2\epsi$
and the  $2n(n-1)$ short roots
$\pm\epsi \pm\epsj , ~~(i< j)$, totalling $2n^2.$
 For the simple system $\Pi$ take
$$
\alpha_1=\epsilon_1-\epsilon_2,
\alpha_2=\epsilon_2-\epsilon_3,\ldots,
\alpha_{n-1}=\eps_{n-1}-\epsn,
\alpha_n=2\epsn .
$$
Then the highest root
$$
\tilde{\alpha}=2\epsilon_1 .
$$
 The Weyl group $W$ is the semidirect product of $S_n$ (which permutes
$\epsi$) and $(\Z/2)^n$ (acting by sign changes on the $\epsi$), the latter
normal in $W$.

 The fundamental coweights are
$$
\boiv=\epsilon_1+\epsilon_2+\ldots+\epsi ,
~~
1\leq i<n~ .
$$
 Now one can repeat word by word the case of $B_n$.
 This completes the proof of part (a).

{\bf Proof of (b)}\qua
 Let now $\Phi$ be of the type $D_n,n\geq 4.$
 Denote by $\epsilon_1,\ldots,\epsilon_n$ the standard basis of $\Rn$.
 Define $\Phi$ to be the set of all vectors  of squared length  2 in the
standard lattice
$\Z\epsilon_1 +\ldots +\Z\epsilon_n$.
 So $\Phi$ consists of the $2n(n-1)$  roots
$\pm\epsi\pm\epsj (1\leq i<j\leq n).$
 For the simple system $\Pi$ take
$$
\alpha_1=\epsilon_1-\epsilon_2,
\alpha_2=\epsilon_2-\epsilon_3,\ldots,
\alpha_{n-1}=\eps_{n-1}-\epsn,
\alpha_n=\eps_{n-1}+\epsn .
$$
Then the highest root
$$
\tilde{\alpha}=\epsilon_1+\epsilon_2 .
$$
 The Weyl group $W$ is the semidirect product of $S_n$ (which permutes
$\epsi$) and $(\Z/2)^{n-1}$ (acting by an even number of sign changes on the
$\epsi$), the latter
normal in $W$.

 The fundamental coweights are
$$
\boiv=\epsilon_1+\epsilon_2+\ldots+\epsi ,
~~
1\leq i<n-2~ ,
$$
$$
\bov_{n-1}=\frac{1}{2}(\epsilon_1+\epsilon_2+\ldots+\epsilon_{n-2}+
\epsilon_{n-1}-\epsilon_n ),
$$
$$
\bov_{n}=\frac{1}{2}(\epsilon_1+\epsilon_2+\ldots+\epsilon_{n-2}+
\epsilon_{n-1}+\epsilon_n ).
$$
Now take the vector 
$$
\omega=\frac{1}{2}(\epsilon_1+\epsilon_2+\ldots+\epsilon_{n-2}
-\epsilon_{n-1}-\epsilon_n ),
$$
which is a coweight since $\omega=u\bonv$ where $u\in W$ acts by 
signs changes on 
$\eps_{n-1},\eps_n.$
Take $j=n-1$, then $\bov_{n-1} \omega=\bov_{n-1}\bonv$
and  it is impossible to replace $u$ by some $w\in W$ so that
$\omega=w\bonv$  and $w$ fixes $\boiv,i< n.$
 Indeed, let $\omega=w\bonv$ and let $w $ be the  identity on 
$\bov_1,\ldots,\bov_{n-1}.$
 Then $w$ fixes $\eps_{n-1}-\eps_n$, from which it follows that 
either $w=1$ or $w$  changes the signs both of $\eps_{n-1},\eps_n.$
Contradiction.
\qed

\eject
\section{Automatic structure for groups acting on Euclid\-ean
buildings  of type $A_n,B_n,C_n$}

\label{main}
\pu
\theorem
\label{main theorem}
{\sl 
{\rm (1)}\qua Let $\Delta$ be any  Euclidean building of one of the
types $A_n,B_n,C_n$, ordered in a standard way  
(see \S\ref{standard ordering} for a definition).
 Then any group acting   freely and cocompactly
on $\Delta$ by type  preserving automorphisms
admits a  biautomatic structure.

{\rm (2)}\qua
 If $\Delta$ is any  Euclidean building of one of the
types $A_n,B_n,C_n$, then any group acting   freely and 
cocompactly on $\Delta$ possesses a  finite index subgroup 
which admits a  biautomatic structure.}\medskip

 Let $G$ be a group satisfying the conditions of the theorem.
 We shall proceed in several steps.
 Firstly we recall the definitions related to an automatic
group theory,
then we establish an isomorphism between the complex $\deltaonespec$ 
and  the Cayley graph of a fundamental groupoid 
$\calg=\pi_1 \left(G\backslash \Delta,G\backslash \deltazerospec 
\right).$
 Making use the combing $\calcspec$ and its properties we prove
the biautomaticity of the groupoid above.
 Finally, we apply a result from \cite{echlpt}, 
\cite{nibloreevescubeaut}, asserting   
that any automorphism group of a biautomatic  groupoid $\calg$
(which is isomorphic to $G$) is biautomatic.


 \pu
{\bf Automatic structures on groups and groupoids}\par
\label{groupoids}
 We shall  use the groupoids technique and since the groups are 
a special case of groupoids, we give all the definitions for
groupoids.
 We summarize here  without proofs the  relevant material on groupoids
from
\cite{echlpt}, \cite{nibloreevescubeaut}.
 A {\it groupoid} $\calg $ is a
category such that the morphism set ${\rm Hom}_\calg (v,w) $ is nonempty
for any  two objects (=vertices) and such that each morphism is
invertible.
 In particular for  any $v\in {\rm Ob}~\calg$ the morphism set
$\calg _v={\rm Hom}_\calg(v,v)$ is a group
and   any group $G$ can be considered as a groupoid with one object,
whose automorphism group is $G$.
 The group $\calg _v$ does not depend on $v$, up to an isomorphism.
 A groupoid is said to be {\it generated }by a set $A$ of
morphisms, if every morphism is a composite of morphisms in $A\cup A^{-1}.
$
Fix a vertex $v_0\in \calg $ as a base point and assume $A$ is a
generating set of morphisms with $A=A^{-1}$ .
 The {\it Cayley graph}
$C\calg =C(\calg ,A,v_0)$ of $\calg $ with respect to a
base point $v_0$ and generating set $A$ is the directed graph with
vertices corresponding to morphisms in $\calg $ with domain (=source)
$v_0$, that is
$$
{\rm Vert}~ C\calg={\rm Hom}(v_0,*)=
\cup\{{\rm Hom}(v_0,v)|v\in{\rm Ob}~\calg\}.
$$
 There
is a directed edge $f\stackrel{a}{\rightarrow }fa$ from the morphism $f$ to
the morphism $fa,\;a\in A$ whenever $fa$ is defined (we write the morphisms
on the right); we give this edge a label $a$.
 In the case of a group the vertex set of a Cayley graph is just a group
itself and we have a usual notion of  a Cayley graph of a group.
 In particular  each edge path in $C\calg$ spells out a word in
$A^{*}$, where as usual, $A^{*}$ denotes the free monoid on $A.$
 And vice versa, for any word $w\in A^*$ and any vertex $f$ there is
a unique path in $C\calg$, beginning in $f$ and having $w$ as its label.
 We put a path-metric  on $C\calg$ by deeming every edge to have unit
length.
 Note that the group $\calg_{v_0}$ acts on $C\calg$ by left
translations.
 The {\it generating graph} $G\calg =G\calg _A$ is a
graph with the same set of vertices as $\calg $ and with edges
corresponding to morphisms in $A$ and labeled by them.

There is a natural projection $p $ from $C\calg $ to
$G\calg$, defined as follows.
 Let the morphism $f:v_0\rightarrow w$
be a vertex of a Cayley graph $C\calg $ then $p (f:v_0\rightarrow w)=w$.
 If $f\stackrel{a}{\rightarrow }fa$
is an edge of $C\calg $ and $f:v_0\rightarrow
w,\;a:w\rightarrow u$ then $p $ sends it to the edge $a:w\rightarrow u$ of
$G\calg $.
 The group $\calg _{v_0}$ gives the group of
covering transformations for $p $ and
$\calg _{v_0}\backslash C(\calg ,A,v)\simeq G\calg .$
 Indeed, the
isomorphism is induced by $p $ and if $f\stackrel{a}{\rightarrow }fa$,
$f_1\stackrel{a_1}{\rightarrow }f_1a_1$ are in the fiber, then
$g=f_1f^{-1}\in \calg _{v_0}.$\ppar

{\bf Automatic structures on groupoids}

 {\ Let $\calg $ be a finitely generated groupoid and $A$ a finite set and
$a\mapsto {\bar a}$ is a map of $A$ to a monoid generating set
${\bar A}\subset \calg $.
 A {\it normal form} for $\calg $ is a subset $L$ of $A^{*}$ satisfying the
following }

(i)\qua $L$ consists of words labelling the paths in $C(\calg,A,v_0)\; $
(that is only composable strings of morphisms are considered,
starting at the base point ${\rm id}\in{\rm Hom}(v_0,v_0)$)

(ii)\qua The natural map $L\rightarrow {\rm Hom}(v_0,*)$ which takes
the word $w=a_1 a_2\cdots a_n$ to the morphism
${\bar a}_1{\bar a}_2\cdots{\bar a}_n\in ~{\rm Hom}(v_0,*)$ is onto.

A {\it rational structure} is a normal form that is a regular language
ie, the set of accepted words for some finite state automaton.
 Recall that a {\it finite state automaton} $\cal M $   with alphabet $A$ is 
a finite directed graph on a vertex set $S$ (called the set of {\it states})
with each edge labeled by an element of $S$ (maybe empty).
 Moreover, a subset of {\it start states} $S_{0}\subset S$ and a subset of
{\it accepted states} $ S_{1}\subset S$ are given. By definition, a word $w$
in the alphabet $A$ is in the
{\it language} $L$ {\it accepted  by} $\cal M$ iff it defines a path
starting from $S_0$ and ending in an accepted state in this graph.
A language is {\it regular} if it is accepted by some finite state automaton.

 We will say that a normal form $L$ has the
``{\it fellow  traveller  property}" if there is a constant $k$ such that
given any normal form words
$v,w\in L$ labelling the paths $\alpha_v,\alpha_w $ in Cayley graph 
$C(\calg ,A,v_0)$ which begin and end at a distance at most one apart,
the distance $d(w(t),v(t)),t=0,1,\dots$ never exceeds $k$.
 A {\it biautomatic structure} for a groupoid $\calg$ is a regular normal form with
the fellow traveller property.

\pu
\label{groupoid is automatic iff vertex group is such one}
\theorem
{\rm (\cite{echlpt}, 13.1.5, \cite{nibloreevescubeaut},
4.1, 4.2)}\qua
{\sl Let }$\calg ${\sl \ be a groupoid and } $v_0$
{\sl  an arbitrary vertex of } $\calg $.
{\sl Then} $\calg $ {\sl admits a  biautomatic 
structure if and only if the  automorphism group }
$\calg _{v_0}$ {\sl \ of } $v_0$ {\sl admits such a structure. }

\pu
{\bf  Groupoid
$\pi_1 \left(G \backslash \Delta,G\backslash \deltazerospec\right)$}\par
\label{fundamental groupoid}
 In this section $G$ is a group acting freely
and cocompactly
on a 
Euclidean building $\Delta$ of of the following types 
$A_n,B_n,C_n$ by automorphisms preserving the standard ordering.

\medskip
{\bf Fundamental groupoid}

The prime  example of a groupoid will be the {\it fundamental groupoid}
$\pi_1 \left( X\right) $
of the path-connected topological space $X$ .
 The set of objects(=vertices)
of $\pi_1\left( X\right) $ is the set of points of $X$
and the morphisms from $x$ to $y$ are homotopy classes of paths beginning
at  $x$ and ending at $y$ .
 The multiplication in $\pi_1\left( X\right) $ is induced by compositions
of paths.
  Given a subset $Y\subset X $ we obtain a
subgroupoid $\pi_1 \left( X,Y\right) $ whose vertices are the points of $Y$
and the morphisms are the same as before.  In particular if $Y$ consists of
a
single point then we get the fundamental group of $X$ based at that point.

\medskip
 {\bf  Generating set of groupoid
$\pi_1 \left( G \backslash \Delta,G\backslash \deltazerospec\right)$}

\pu
 \lemma
\label{generating set for groupoid}
{\sl  Let $A$ be the set of
homotopy classes of the images in $G \backslash \Delta$ of
directed edges of  the graph $\deltaonespec$.
 Then $A$ is a finite set, generating  groupoid
$\calg=\pi_1 \left( G\backslash \Delta,G\backslash \deltazerospec 
\right).$
}

\proof
 This set is finite since by condition $\Delta$ has only finitely
many cells under the action of $G.$
 To prove that $A$ generates $\calg$ take a path $\alpha $ from
$G v$ to $G v',$ then since $\Delta $ is contractible
, see \cite{brownbuildings} and
$G $ acts freely on $\Delta $ the projection
$\Delta \rightarrow G\backslash \Delta \;$ is a
universal cover and there is a unique lift $\widetilde{\alpha }$ of $\alpha
$
 into $\Delta $ which begins at $v.$ This lift ends at a translate
 $gv'$ of $v'$ where $g$ is determined by the homotopy class of $\alpha .$
 Moreover, any path in $\Delta$ from $v$ to $gv'$ will project to
a
path in $G\backslash \Delta$ which is homotopic to $\alpha .$
 In particular  the path from
$\calc_{spec}$ which crawls from $v$ to $gv'$ is homotopic to 
$\widetilde{\alpha}$.
 Since this path is the edge path in $\deltaonespec$ it projects into
$G\backslash\Delta$ as a composition of
homotopy classes of the images in $G \backslash \Delta $ of
directed edges of  the graph $\deltaonespec,$ which is a product of
 morphisms from $A$.
 This means that $A$ is a generating set for
$\calg$.
\qed

\medskip
{\bf Labelling the  graph $\deltaonespec$}

 Consider $\deltaonespec$ as a directed graph and label the edge $a$
by an element $G a\in A$.

 \pu
\lemma
\label{iso cayley and skeleton}
{\sl 
$\deltaonespec$ as a labeled graph is isomorphic to a Cayley graph
$C\calg$ of a fundamental groupoid ~
$\calg=\pi_1 \left( G\backslash \Delta,G\backslash \deltazerospec 
\right).$
}

\proof
Fix a base vertex $v_0$ in $\deltaonespec$ and consider
$G v_0$ as a base point in $G \backslash \Delta $.
 A vertex in $\calg$ is  a homotopy class $[f]$ of  paths from 
$Gv_0$ to some $G v$.
  There is a   unique lift $\widetilde{f}$ of $f $ into $\Delta $ which
begins at $v_0.$
  We send the vertex $[f]$ to the end of the path $\widetilde{f}$.
  Now if $[f]\stackrel{[a]}{\rightarrow }[f][a]$ is an edge in 
$C\calg$ then there are unique lifts $\widetilde{f},\widetilde{a}$ of 
$[f],[a]$ to $\Delta $ such that $\widetilde{a}$ starts at the end of 
$\widetilde{f}$.
 We map the edge $[f]\stackrel{[a]}{\rightarrow }[f][a]$
to the edge $\widetilde{a}$, labeled by $G \widetilde{a}$.
 This defines an isomorphism as required.
\qed

 \pu
{\bf Language}\par
\label{language}
Recall that we label  directed edges in $\deltaonespec$
in a $G$--equivariant way  by $A$,
so each path from $\calcspec$ spells out a word in
$A^{*}.$
 Define a language $L$ to be the subset of $A^{*}$ which is given by
all words which label the paths from combing $\calcspec$ starting at the
basepoint $v_0$.
 It follows from the above discussion that we have a bijective map from $L$ to
morphisms in 
$\calg=\pi_1 \left( G\backslash \Delta,G\backslash \deltazerospec 
\right).$

\medskip
\lemma
\label{L is regular}
{\sl The language }$L$  {\sl over } $A$ {\sl determined by the combing }
$\calc_{spec}$ {\sl is regular.}

\proof
(cf \cite{nibloreevescubeaut}, 6.1)\qua
  We shall construct a non-deterministic finite state automaton
$M$ over $A$ which has $L$ as the set of acceptable words.
 The set of states of $M$ is $A;$ all states are initial states and all
states are acceptable states.
 There is a transition labelled by $G e_1$ from $G e_1$ to
$G e_2$ if there are ordered  
edges $e_1',e_2'\in \deltaonespec$ in
$G e_1,G e_2$ respectively, such that $e_2'$
starts at the tail of $e_1'$ and
one of the following  conditions holds:

1)\qua both $e_1',e_2'$ constitute a geodesic
linear path of the length 2
in $\deltaonespec$, that is 
a local geodesicity condition is
satisfied in the common vertex, see \ref{glg}.

2)\qua If $e_1'$ is of type $i$ and $e_2'$ is of type $j$, and if 
$e_1''$ is the last edge in $\deltaone$, lying on the segment
$e_1'$ and $e_2''$ is the first edge in $\deltaone$ lying on
$e_2'$,
then
$(e_1'',e_2'')=(\frac{\boiv}{c_i},\frac{\bojv}{c_j}).$

Since the condition defining the transitions is the same as in
the local
description of $\calc$ in \S \ref{recursive}, the language $L$ is 
precisely the language
accepted by the finite state automaton.

\pu
{\bf  Finishing the proof of the theorem \ref{main theorem}}\par
\label{proofmain}
 By Theorem \ref{groupoid is automatic iff vertex group is such one}
it is enough to show that the fundamental groupoid
$\calg=$\break$\pi_1(G\backslash \Delta,G\backslash
\deltazerospec)$
is automatic.
  Fix a base vertex  $G v_0$ for $G\backslash \Delta$,
where $v_0\in \deltazerospec$.
  Let $A$ be an
 alphabet which is in one one correspondence with a finite generating set
$G\backslash \deltaonespec$ of $\calg$, see
\ref{generating set for groupoid}.
  Let $L\subseteq A^{*}$ be a language consisting of words which are
spelled out from  paths of $\calg$.
  By the construction of a combing, see \ref{def combing},
$L$ surjects onto $\calg(v_0,*)$. 
 By \S\ref{L is regular}, $L$ is regular. 
 By \S\ref{ftp} it  satisfies $k$--fellow traveller property.
 Hence, in view of an isomorphism of $A$--labelled  graphs
$\deltaonespec\simeq C\calg$ we get that $\calg $  is biautomatic.

 To prove the second assertion 
just note that the set of all orderings of a Euclidean building is finite 
and any group acting simplicially on a building, acts also
on this finite set of orderings, hence it contains a subgroup
of finite index which preserves any ordering,
that is acts by a type preserving automorphisms.
\qed

\pu
\remarks
 Actually, one can derive easily from 
Sections \ref{qi1},\ref{qi2},\ref{combing is quasigeodesic}
that the structures we have built are quasigeodesic ones.
 On the other hand, as Prof W Neumann has pointed out to us, 
every  (synchronous) automatic structure contains
a sublanguage which is an  (synchronous) automatic structure with 
the uniqueness property (\cite{echlpt} 2.5.1) and  it follows  
(\cite{echlpt} 3.3.4) that (synchronous) automatic structures with 
uniqueness are always quasigeodesic.

\end{document}